\documentclass[a4paper,12pt,leqno]{article}
\usepackage{longtable}
\usepackage{latexsym}

\title{The classification of smooth toric weakened Fano 3-folds}

\date{}

\author{{\sc Hiroshi Sato}\thanks{Partly supported by the Grant-in-Aid for JSPS Fellows, The Ministry of Education, Science, Sports and Culture, Japan.
\newline
\hspace*{1.5em} {\em $2000$ Mathematics Subject Classification\/}.
Primary 14M25;
Secondary 14J30, 14J45.}}

\newtheorem{Thm}{{Theorem}}[section]
\newtheorem{Prop}[Thm]{{Proposition}}
\newtheorem{Cor}[Thm]{{Corollary}}

\newtheorem{Lem}[Thm]{{Lemma}}

\newtheorem{Def}[Thm]{{Definition}}
\newtheorem{Rem}[Thm]{{Remark}}

\newcommand{\proof}{Proof. \quad}
\newcommand{\qed}{\hfill $\Box$}

\newcommand{\G}{\mathop{\rm G}\nolimits}

\newcommand{\PC}{\mathop{\rm PC}\nolimits}

\begin{document}
\maketitle

\begin{abstract}                                                  

We completely classify toric weakened Fano $3$-folds, that is, smooth toric weak Fano $3$-folds which are not Fano but are deformed to smooth Fano $3$-folds. There exist exactly $15$ toric weakened Fano $3$-folds up to isomorphisms.

\end{abstract}

\section{Introduction}\label{introduction}    

\thispagestyle{empty}

\hspace{5mm} A {\em Fano} (resp. {\em weak Fano}) variety $X$ is a smooth projective variety whose anti-canonical divisor $-K_X$ is ample (resp. nef and big). Minagawa \cite{minagawa1} introduce the concept of {\em weakened Fano} variety in connection with ``Reid's Fantasy'' for weak Fano $3$-folds. A weak Fano variety $X$ is called a weakened Fano variety if it is not Fano but is deformed to Fano under a small deformation (see Definition \ref{defweakened}). In this paper, we consider the classification problem of weakened Fano $3$-folds for the case of toric varieties. As a result, we can determine the structures of toric weakened Fano $3$-folds using a result of Minagawa \cite{minagawa1}, \cite{minagawa2}. There exist exactly $15$ toric weakened Fano $3$-folds up to isomorphisms (see Theorem \ref{class}). There are three cases: $(1)$ ${\bf P}^1\times X'$, where $X'$ is a toric weak del Pezzo surface but not a del Pezzo surface, $(2)$ toric del Pezzo surface bundles over ${\bf P}^1$ and $(3)$ toric weakened del Pezzo surface bundles over ${\bf P}^1$.

The content of this paper is as follows: Section \ref{primsection} is a section for preparation. We review the basic concepts such as the toric Mori theory, primitive collections and primitive relations. In Section \ref{surface}, we give the classification of toric weak del Pezzo surfaces. This is necessary for the classification of smooth toric weakened Fano $3$-folds. In Section \ref{mainsec}, we consider the classification of smooth toric weakened Fano $3$-folds.

The author wishes to thank Professor Tatsuhiro Minagawa for advice and encouragement. The problem of the classification of toric weakened Fano $3$-folds was posed by Professor Shihoko Ishii. The author also thanks her for advice.

\section{Primitive collections and primitive relations}\label{primsection}

\hspace{5mm} In this section, we review the concepts of smooth complete toric varieties using the notion of primitive collections and primitive relations. See Batyrev \cite{batyrev3}, \cite{batyrev4} and Sato \cite{sato1} more precisely. For fundamental properties of the toric geometry, see Fulton \cite{fulton1} and Oda \cite{oda2}. We work over the complex field {\bf C} throughout this paper.

\begin{Def}

{\rm
Let $X$ be a smooth complete toric $d$-fold, $\Sigma$ the corresponding fan in $N:={\bf Z}^d$ and $\G(\Sigma)\subset N$ the set of primitive generators of $1$-dimensional cones in $\Sigma$. A subset $P\subset\G(\Sigma)$ is called a {\em primitive collection} of $\Sigma$ if $P$ does not generate a cone in $\Sigma$ while any proper subset of $P$ generates a cone in $\Sigma$. We denote by $\PC(\Sigma)$ the set of primitive collections of $\Sigma$.
}

\end{Def}

Let $P=\{x_1,\ldots,x_m\}$ be a primitive collection of $\Sigma$. Then there exists a unique cone $\sigma(P)$ in $\Sigma$ such that $x_1+\cdots+x_m$ is contained in the relative interior of $\sigma(P)$, because $X$ is complete. So we get an equality
$$x_1+\cdots+x_m=a_1y_1+\cdots+a_ny_n,$$
where $y_1,\ldots,y_n$ are the generators of $\sigma(P)$, that is, $\sigma(P)\cap\G(\Sigma)$ and $a_1,\ldots,a_n$ are positive integers. We call this equality the {\em primitive relation} of $P$. Thus we get an element $r(P)$ in $A_{1}(X)$ for any primitive collection $P\in\PC(\Sigma)$, where $A_{1}(X)$ is the group of $1$-cycles on $X$ modulo rational equivalences. We define the {\em degree} of $P$ as $\deg P:=\left( -K_X\cdot r(P)\right)=m-(a_1+\cdots+a_n)$. The following is important.

\begin{Prop}[Batyrev \cite{batyrev3}, Reid \cite{reid1}]\label{toriccone}

Let $X$ be a smooth projective toric variety and $\Sigma$ the corresponding fan. Then
$${\bf NE}(X)=\sum_{P\in\PC(\Sigma)}{\bf R}_{\geq 0}r(P),$$
where ${\bf NE}(X)$ is the Mori cone of $X$.

\end{Prop}

A primitive collection $P$ is called an {\em extremal} primitive collection when $r(P)$ is contained in an extremal ray of ${\bf NE}(X)$.

\section{Toric weak Fano varieties}\label{surface}

\hspace{5mm} In this section, we review the concepts of toric Fano varieties and toric weak Fano varieties. Especially, we give the classification of toric weak del Pezzo surfaces.

\begin{Def}

{\rm
Let $X$ be a smooth projective algebraic variety. Then $X$ is called a {\em Fano} variety (resp. {\em weak Fano} variety), if its anti-canonical divisor $-K_X$ is ample (resp. nef and big).
}

\end{Def}

For the toric case, Fano varieties and weak Fano varieties are characterized as follows.

\begin{Prop}[Batyrev \cite{batyrev4}, Sato \cite{sato1}]\label{odafano}

Let $X$ be a smooth projective toric variety and $\Sigma$ the corresponding fan. Then $X$ is Fano $($resp. weak Fano$)$ if and only if $\deg P>0$ $($resp. $\deg P\geq 0)$ for any primitive collection $P\in\PC(\Sigma)$.

\end{Prop}

By Propositions \ref{toriccone} and \ref{odafano}, the following holds. This is very useful in Section \ref{mainsec}.

\begin{Cor}\label{primfano}

Let $X$ be a toric weak Fano variety and $\Sigma$ the corresponding fan. Then for any primitive collection $P=\{x_1,x_2\}\in\PC(\Sigma)$, the corresponding primitive relation is one of the following.
\begin{enumerate}
\item $x_1+x_2=0$.
\item $x_1+x_2=ay$ $(y\in\G(\Sigma)$, and either $a=1$ or $2)$.
\item $x_1+x_2=y_1+y_2$ $\left( y_1,\ y_2\in\G(\Sigma)\right)$.
\end{enumerate}

\end{Cor}

There exists a one-to-one correspondence between smooth toric weak del Pezzo surfaces, that is, $2$-dimensional smooth toric weak Fano varieties and Gorenstein toric del Pezzo surfaces (see Section 6 in Sato \cite{sato1}). Since Gorenstein toric del Pezzo surfaces are classified (see Koelman \cite{koelman1}), we can completely classify smooth toric weak del Pezzo surfaces. We give all the smooth toric weak del Pezzo surfaces by giving the elements of $\G(\Sigma)$ (see Table 1). This classification is necessary in the following section.
\setlongtables
\bigskip

\bigskip

\newpage

\begin{center}
Table 1: Smooth toric weak del Pezzo surfaces
\end{center}

\begin{longtable}{|c|c|c|}

\hline
\endhead
\hline
\endfoot

 &  $\G(\Sigma)$ & notation \\ \hline
(1) & $\pmatrix{1\cr 0},\ \pmatrix{0\cr 1},\ \pmatrix{-1\cr -1}.$ & ${\bf P}^{2}$ \\ \hline
(2) & $\pm\pmatrix{1\cr 0},\ \pm\pmatrix{0\cr 1}.$ & ${\bf P}^{1}\times{\bf P}^1$ \\ \hline
(3) & $\pm\pmatrix{1\cr 0},\ \pmatrix{0\cr 1},\ \pmatrix{1\cr -1}.$ & $F_1$ \\ \hline
(4) & $\pm\pmatrix{1\cr 0},\ \pmatrix{1\cr \pm1}.$ & $F_2$ \\ \hline
(5) & $\pm\pmatrix{1\cr 0},\ \pm\pmatrix{0\cr 1},\ \pmatrix{1\cr 1}.$ & $S_7$ \\ \hline
(6) & $\pm\pmatrix{1\cr 0},\ \pmatrix{1\cr \pm 1},\ \pmatrix{0\cr 1}.$ & $W_3$ \\ \hline
(7) & $\pm\pmatrix{1\cr 0},\ \pm\pmatrix{0\cr 1},\ \pm\pmatrix{1\cr 1}.$ & $S_6$ \\ \hline
(8) & $\pm\pmatrix{1\cr 0},\ \pm\pmatrix{0\cr 1},\ \pmatrix{1\cr \pm 1}.$ & $W^1_4$ \\ \hline
(9) & $\pm\pmatrix{1\cr 0},\ \pmatrix{1\cr \pm1},\ \pmatrix{0\cr 1},\ \pmatrix{-1\cr 1}.$ & $W^2_4$ \\ \hline
(10) & $\pm\pmatrix{1\cr 0},\ \pmatrix{1\cr \pm1},\ \pmatrix{0\cr 1},\ \pmatrix{1\cr 2}.$ & $W^3_4$ \\ \hline
(11) & $\pm\pmatrix{1\cr 0},\ \pm\pmatrix{0\cr 1},\ \pmatrix{1\cr \pm 1},\ \pmatrix{-1\cr 1}.$ & $W^1_5$ \\ \hline
(12) & $\pm\pmatrix{1\cr 0},\ \pm\pmatrix{0\cr 1},\ \pmatrix{1\cr \pm 1},\ \pmatrix{1\cr -2}.$ & $W^2_5$ \\ \hline
(13) & $\pm\pmatrix{1\cr 0},\ \pm\pmatrix{0\cr 1},\ \pmatrix{1\cr \pm 1},\ \pmatrix{-1\cr \pm 1}.$ & $W^1_6$ \\ \hline
(14) & $\pm\pmatrix{1\cr 0},\ \pm\pmatrix{0\cr 1},\ \pmatrix{1\cr \pm 1},\ \pmatrix{1\cr -2},\ \pmatrix{-1\cr 1}.$ & $W^2_6$ \\ \hline
(15) & $\pm\pmatrix{1\cr 0},\ \pm\pmatrix{0\cr 1},\ \pmatrix{1\cr \pm 1},\ \pmatrix{1\cr \pm 2}.$ & $W^3_6$ \\ \hline
(16) & $\pm\pmatrix{1\cr 0},\ \pm\pmatrix{0\cr 1},\ \pmatrix{1\cr \pm 1},\ \pmatrix{1\cr -2},\ \pmatrix{-1\cr 1},\ \pmatrix{-2\cr 1}.$ & $W_7$ \\ \hline

\end{longtable}

\bigskip

In Table 1, we denote by $F_a$ the Hirzebruch surface of degree $a$, while we denote by $S_n$ the del Pezzo surface of degree $n$.

\section{Classification}\label{mainsec}

\hspace{5mm} Minagawa defined the following concept in connection with ``Reid's fantasy'' for weak Fano $3$-folds.

\begin{Def}[Minagawa \cite{minagawa1}]\label{defweakened}

{\rm
Let $X$ be a weak Fano variety. Then $X$ is called a {\em weakened Fano} variety if
\begin{enumerate}
\item $X$ is not a Fano variety and
\item there exists a small deformation $\varphi :\mathcal{X}\rightarrow\Delta_{\epsilon}:=\left\{t\in{\bf C}\,\left|\,|t|<\epsilon<<1\right.\right\}$ such that $\mathcal{X}_{0}:=\varphi^{-1}(0)\cong X$ while $\mathcal{X}_t:=\varphi^{-1}(t)$ is a Fano variety for any $t\in\Delta_{\epsilon}\setminus\{0\}$. 
\end{enumerate}

}

\end{Def}

\begin{Rem}

{\rm
Let $X$ be a weak del Pezzo surface. If $X$ is not a del Pezzo surface then $X$ is a weakened del Pezzo surface.
}

\end{Rem}

\begin{Rem}

{\rm
Weakened Fano $3$-folds of Picard number two are studied in Minagawa \cite{minagawa1}.
}

\end{Rem}

The main purpose of this paper is to classify toric weakened Fano $3$-folds.

Minagawa characterized weakened Fano $3$-folds using the notion of primitive contractions.

\begin{Thm}[Minagawa \cite{minagawa1}, \cite{minagawa2}]\label{charweakened}

Let $X$ be a weak Fano $3$-fold and not a Fano $3$-fold. Then $X$ is a weakened Fano $3$-fold if and only if every primitive crepant contraction $f:X\rightarrow\overline{X}$ is a divisorial contraction which contracts a divisor $E\subset X$ to a curve $\overline{C}\subset\overline{X}$ such that
\begin{enumerate}
\item $f|_{E}:E\rightarrow \overline{C}$ is a ${\bf P}^1$-bundle structure,
\item $\overline{C}\cong{\bf P}^1$ and
\item $(-K_{\overline{X}}\cdot \overline{C})=2$.
\end{enumerate}
Such contractions are called {\em $(0,2)$-type} contractions.

\end{Thm}

We study contractions of $(0,2)$-type for the case of toric varieties. Let $X$ be a toric weakened Fano $3$-fold and $\Sigma$ the corresponding fan in $N\cong{\bf Z}^3$. By the assumptions $f$ is crepant and $f$ contracts a divisor to a curve, without loss of generalities, we can assume $\Sigma$ contains four $3$-dimensional cones
$$\sigma_1={\bf R}_{\geq 0}x_0+{\bf R}_{\geq 0}x_++{\bf R}_{\geq 0}y_+,\ \sigma_2={\bf R}_{\geq 0}x_0+{\bf R}_{\geq 0}x_++{\bf R}_{\geq 0}y_-,$$
$$\sigma_3={\bf R}_{\geq 0}x_0+{\bf R}_{\geq 0}x_-+{\bf R}_{\geq 0}y_+\mbox{ and }\sigma_4={\bf R}_{\geq 0}x_0+{\bf R}_{\geq 0}x_-+{\bf R}_{\geq 0}y_-,$$
where
$$x_0=\pmatrix{1\cr 0\cr 0},\ x_+=\pmatrix{1\cr 1\cr 0},\ x_-=\pmatrix{1\cr -1\cr 0},\ y_+=\pmatrix{0\cr 0\cr 1}\mbox{ and }y_-=\pmatrix{\alpha\cr\beta\cr -1}$$for some integers $\alpha$ and $\beta$. In this case, $E$ is the toric divisor corresponding to the $1$-dimensional cone ${\bf R}_{\geq 0}x_0$ in $\Sigma$. The conditions $(1)$ and $(2)$ are automatically satisfied. On the other hand, the set of maximal cones of the fan $\overline{\Sigma}$ corresponding to $\overline{X}$ is
$$\left(\left\{\mbox{the maximal cones of }\Sigma\right\}\setminus\left\{\sigma_1,\sigma_2,\sigma_3,\sigma_4\right\}\right)\cup\left\{\overline{\sigma_1}:=\sigma_1\cup\sigma_3,\ \overline{\sigma_2}:=\sigma_2\cup\sigma_4\right\}.$$
$\overline{C}$ is the torus invariant curve corresponding to the $2$-dimensional cone ${\bf R}_{\geq 0}x_++{\bf R}_{\geq 0}x_-$ in $\overline{\Sigma}$. Put $C$ be the torus invariant curve corresponding to the cone ${\bf R}_{\geq 0}x_0+{\bf R}_{\geq 0}x_+$ on $X$. Then we have $\left(-K_{\overline{X}}\cdot \overline{C}\right)=\left(-K_{\overline{X}}\cdot f_{\ast}C\right)=\left(f^{\ast}(-K_{\overline{X}})\cdot C\right)=\left(-K_X\cdot C\right)$, because $\overline{C}=f_{\ast}C$ in $A_{1}(\overline{X})$ while $f^{\ast}K_{\overline{X}}=K_X$ in ${\rm Pic}(X)$. Let $D_0$, $D_+$, $D_-$, $H_+$ and $H_-$ be the toric divisors on $X$ corresponding to $x_0$, $x_+$, $x_-$, $y_+$ and $y_-$, respectively. Then the equalities
$$D_0+D_++D_-+\alpha H_-+\mathcal{D}_1=0,\ D_+-D_-+\beta H_-+\mathcal{D}_2=0\mbox{ and}$$
$$H_+-H_-+\mathcal{D}_3=0$$
hold in ${\rm Pic}(X)$, where $\mathcal{D}_1$, $\mathcal{D}_2$ and $\mathcal{D}_3$ are linear combinations of prime toric devisors on $X$ other than $D_0$, $D_+$, $D_-$, $H_+$ and $H_-$. Note $(H_+\cdot C)=(H_-\cdot C)=1$ and $(D_-\cdot C)=(\mathcal{D}_1\cdot C)=(\mathcal{D}_2\cdot C)=(\mathcal{D}_3\cdot C)=0$. By these equalities, we have $(D_+\cdot C)=-\beta$ and $(D_0\cdot C)=\beta-\alpha$. Thus we have $\left(-K_{\overline{X}}\cdot \overline{C}\right)=\left(-K_X\cdot C\right)=(D_0\cdot C)+(D_+\cdot C)+(D_-\cdot C)+(H_+\cdot C)+(H_-\cdot C)=(\beta-\alpha)-\beta+2=2-\alpha.$ So by the assumption $\left(-K_{\overline{X}}\cdot \overline{C}\right)=2$, we have $\alpha=0$ and $E$ is isomorphic to the Hirzebruch surface $F_{\beta}$ of degree $\beta$. Moreover, since $X$ is a weak Fano $3$-fold, we have $-2\leq\beta\leq 2$. We may assume $0\leq\beta\leq 2$. As a result, we get the following.

\begin{Prop}\label{zerotwo}

Let $X$ be a toric weak Fano $3$-fold and $\Sigma$ the corresponding fan. There exists a $(0,2)$-type contraction $f:X\rightarrow\overline{X}$ if and only if, by some automorphism of $N={\bf Z}^3$, $\Sigma$ contains four $3$-dimensional cones
$$\sigma_1={\bf R}_{\geq 0}x_0+{\bf R}_{\geq 0}x_++{\bf R}_{\geq 0}y_+,\ \sigma_2={\bf R}_{\geq 0}x_0+{\bf R}_{\geq 0}x_++{\bf R}_{\geq 0}y_-,$$
$$\sigma_3={\bf R}_{\geq 0}x_0+{\bf R}_{\geq 0}x_-+{\bf R}_{\geq 0}y_+\mbox{ and }\sigma_4={\bf R}_{\geq 0}x_0+{\bf R}_{\geq 0}x_-+{\bf R}_{\geq 0}y_-,$$
where
$$x_0=\pmatrix{1\cr 0\cr 0},\ x_+=\pmatrix{1\cr 1\cr 0},\ x_-=\pmatrix{1\cr -1\cr 0},\ y_+=\pmatrix{0\cr 0\cr 1}\mbox{ and }y_-=\pmatrix{0\cr a\cr -1},$$
and $0\leq a\leq 2$. Especially, the exceptional divisor $E$ of $f$ corresponds to ${\bf R}_{\geq 0}x_0\in\Sigma$ and $E\cong F_a$.

\end{Prop}

The following is fundamental to consider the classification of toric weakened Fano $3$-folds.

\begin{Lem}\label{key}

Let $X$ be a toric weakened Fano $3$-fold and $\Sigma$ the corresponding fan. Then for any primitive collection $P=\{x_1,x_2\}\in\PC(\Sigma)$ whose primitive relation is $x_1+x_2=y_1+y_2$, where $\sigma(P)\cap\G(\Sigma)=\{y_1,y_2\}$, we have ${\bf R}_{\geq 0}x_1+{\bf R}_{\geq 0}x_2\supset{\bf R}_{\geq 0}y_1+{\bf R}_{\geq 0}y_2$.

\end{Lem}

\proof
Suppose ${x_1,x_2,y_1,y_2}$ generates a $3$-dimensional cone. If $P$ is extremal, then the corresponding primitive contraction is a small contraction. So by Theorem \ref{charweakened}, $P$ is not extremal. Therefore, there exist extremal primitive collections $P_1,\ldots,P_n$ such that
$$r(P)=\sum_{i=1}^{n}a_ir(P_i),$$
where $a_1,\ldots,a_n$ are positive integers and $n\geq 2$. Since $\deg P=(-K_X\cdot r(P))=0$, for any $1\leq i\leq n$, we have $\deg P_i=(-K_X\cdot r(P_i))=0$ and the corresponding primitive crepant contractions are $(0,2)$-type by Theorem \ref{charweakened}. Let $x_i'+x_i''=2y_i'$ be the corresponding primitive relation of $P_i$ for any $1\leq i\leq n$. Then there exist $1\leq j,k\leq n$ such that $y_1=y_j'$ and $y_2=y_k'$. This is impossible because $P_j$ and $P_k$ are extremal. Thus we have ${\bf R}_{\geq 0}x_1+{\bf R}_{\geq 0}x_2\supset{\bf R}_{\geq 0}y_1+{\bf R}_{\geq 0}y_2$.
\qed

\begin{Cor}\label{keycor}

Let $X$ be a toric weakened Fano $3$-fold and $\Sigma$ the corresponding fan. For any primitive collecion $P=\{x_1,x_2\}\in\PC(\Sigma)$ such that $x_1+x_2\neq 0$, there exists an element $z\in\G(\Sigma)$ such that $z$ is contained in the relative interior of ${\bf R}_{\geq 0}x_1+{\bf R}_{\geq 0}x_2$.

\end{Cor}

\proof
This is obvious by Corollary \ref{primfano} and Lemma \ref{key}.
\qed

\bigskip

Now we can start the classification of toric weakened Fano $3$-folds.

Let $X$ be a toric weakened Fano $3$-fold and $\Sigma$ the corresponding fan. Since $X$ is not Fano, there exists a $(0,2)$-type contraction. We use the notation as in Proposition \ref{zerotwo}. Put
$$\mathcal{S}_+=\left\{\left.\alpha x_0+\beta x_+\in N_{\bf R}:=N\otimes_{\bf Z}{\bf R}\cong{\bf R}^3\; \right|\;\alpha,\beta\in{\bf R},\beta\geq 0\right\},$$
$$\mathcal{S}_-=\left\{\left.\alpha x_0+\beta x_-\in N_{\bf R}\; \right|\;\alpha,\beta\in{\bf R},\beta\geq 0\right\},$$
$$\mathcal{T}_+=\left\{\left.\alpha x_0+\beta y_+\in N_{\bf R}\; \right|\;\alpha,\beta\in{\bf R},\beta\geq 0\right\},$$
$$\mathcal{T}_-=\left\{\left.\alpha x_0+\beta y_-\in N_{\bf R}\; \right|\;\alpha,\beta\in{\bf R},\beta\geq 0\right\},$$
$$\mathcal{S}=\mathcal{S}_+\cup \mathcal{S}_-\mbox{ and }\mathcal{T}=\mathcal{T}_+\cup \mathcal{T}_-.$$
$\mathcal{S}_+$, $\mathcal{S}_-$, $\mathcal{T}_+$, $\mathcal{T}_-$, $\mathcal{S}$ and $\mathcal{T}$ are connected subsets in $N_{\bf R}$. The following holds.

\begin{Lem}\label{wataru}

$I:=\G(\Sigma)\setminus\{x_0,x_+,x_-,y_+,y_-\}$ is contained in either $\mathcal{S}$ or $\mathcal{T}$.

\end{Lem}

\proof
First we show $\G(\Sigma)\subset \mathcal{S}\cup \mathcal{T}$. Suppose there exists $z\in I$ such that $z\not\in \mathcal{S}\cup \mathcal{T}$. Since $\{ x_0,z\}$ is a primitive collection, there exists $z'\in\G(\Sigma)$ such that $z'$ is contained in the relative interior of ${\bf R}_{\geq 0}x_0+{\bf R}_{\geq 0}z$ by Corollary \ref{keycor}. We can replace $z$ by $z'$ and do this discussion again. This contradicts the fact $\G(\Sigma)$ is a finite set. So we have $\G(\Sigma)\subset \mathcal{S}\cup \mathcal{T}$.

Next suppose there exist $z_1$ and $z_2\in I$ such that $x_0+z_1\neq 0$, $x_0+z_2\neq 0$, $z_1\in \mathcal{S}_+$ and $z_2\in \mathcal{T}_+$ (the other cases are similar). By the similar discussion as above, we can choose such $z_1$ and $z_2$ as $\{x_+,z_1\}$ and $\{y_+,z_2\}$ generate $2$-dimensional cones in $\Sigma$. On the other hand, $\{z_1,z_2\}$ generates a $2$-dimensional cone in $\Sigma$ by Corollary \ref{keycor} and the above discussion. Therefore, either $\{x_+,z_2\}$ or $\{y_+,z_1\}$ is a primitive collection. By Corollary \ref{keycor}, this contradicts the fact $\G(\Sigma)\subset \mathcal{S}\cup \mathcal{T}$.
\qed

\begin{Lem}\label{split}

If $\{y_+,y_-\}$ $($resp. $\{x_+,x_-\})$ is a primitive collection and $I\subset\mathcal{S}$ $($resp. $I\subset\mathcal{T}$$)$, then for any primitive collection $P\in\PC(\Sigma)$ such that $P\neq\{y_+,y_-\}$ $($resp. $P\neq\{x_+,x_-\})$, we have $P\cap\{y_+,y_-\}=\emptyset$ $($resp. $P\cap\{x_+,x_-\}=\emptyset)$.

\end{Lem}

\proof
We show this lemma for the case $I\subset\mathcal{S}$. Another case is similar.

Suppose $P\in\PC(\Sigma)$ be a primitive collection such that $P\neq\{y_+,y_-\}$. By Corollary \ref{keycor} and Lemma \ref{wataru}, we have $\# P=3$. So, without loss of generalities, we may assume $P=\{y_+,z_1,z_2\}$, where $z_1,z_2\in I$. Since $y_+$, $z_1$ and $z_2$ are linearly independent over ${\bf R}$, there exists a element of $\G(\Sigma)$ in the interior of ${\bf R}_{\geq 0}y_++{\bf R}_{\geq 0}z_1+{\bf R}_{\geq 0}z_2$. This contradicts Lemma \ref{wataru}.
\qed

\begin{Rem}\label{f1only}

{\rm
In Lemma \ref{split}, the condition $\{x_+,x_-\}$ is a primitive collection holds automatically. Moreover, if $E\cong{\bf P}^1\times{\bf P}^1$ or $E\cong F_2$, that is, $a=0$ or $a=2$, then the condition $\{y_+,y_-\}$ is a primitive collection also holds automatically. In these cases, $y_++y_-=0$ and $y_++y_-=2z$ $(z\in\G(\Sigma))$ are the corresponding primitive relations, respectively.
}

\end{Rem}

\begin{Cor}\label{bdl}

Under the assumption in Lemma {\rm\ref{split}}, if the Picard number of $X$ is not three, then one of the following holds.
\begin{enumerate}
\item $I\subset\mathcal{S}$. Moreover, $X$ is a toric surface bundle over ${\bf P}^1$ such that the fan corresponding to a fiber is in $\mathcal{S}$.
\item $I\subset\mathcal{T}$ and $E\cong{\bf P}^1\times{\bf P}^1$. Moreover, $X$ is a toric surface bundle over ${\bf P}^1$ such that tha fan corresponding to a fiber is in $\mathcal{T}$.
\end{enumerate}

\end{Cor}

\proof
See Proposition 4.1, Theorem 4.3 and Corollary 4.4 in Batyrev \cite{batyrev3}.
\qed

\bigskip

By these results, we can complete the classification of toric weakened Fano $3$-folds. We split the classification into three cases, that is, (I) $E\cong{\bf P}^1\times{\bf P}^1$, (II) $E\cong F_1$ and (III) $E\cong F_2$.

\bigskip
$({\rm I})$ $E\cong{\bf P}^1\times{\bf P}^1\ (a=0).$
\bigskip

If $I\subset\mathcal{S}$, then $X$ is isomorphic to ${\bf P}^1\times X'$ by $(1)$ in Corollary \ref{bdl}, where $X'$ is a toric weak del Pezzo surface but not a toric del Pezzo surface. In this case, by the classification of toric weak del Pezzo surfaces in Section \ref{surface}, there exist exactly $11$ toric weakened Fano $3$-folds up to isomorphisms: ${\bf P}^1\times F_2$, ${\bf P}^1\times W_3$, ${\bf P}^1\times W^1_4$, ${\bf P}^1\times W^2_4$, ${\bf P}^1\times W^3_4$, ${\bf P}^1\times W^1_5$, ${\bf P}^1\times W^2_5$, ${\bf P}^1\times W^1_6$, ${\bf P}^1\times W^2_6$, ${\bf P}^1\times W^3_6$ and ${\bf P}^1\times W_7$ (see Table $1$ in Section \ref{surface}).

If $I\subset\mathcal{T}$, then $X$ is a toric surface bundle over ${\bf P}^1$ by $(2)$ in Corollary \ref{bdl}. A fiber $X'$ of this bundle structure corresponds to the $2$-dimensional fan in $\mathcal{T}$ and $X'$ is a toric weak del Pezzo surface. Therefore, by the classification of toric weak del Pezzo surfaces in Section \ref{surface}, in this case, we get two new toric weakened Fano $3$-folds $X_3^0$ and $X_4^0$ whose Picard numbers are $3$ and $4$, respectively. The other cases are impossible because there exists a crepant contraction which is not of $(0,2)$-type. Put
$$z_1=\pmatrix{-1\cr 0\cr 1},\ z_2=\pmatrix{-1\cr 0\cr 0},$$
$\Sigma_3^0$ the fan corresponding to $X_3^0$ and $\Sigma_4^0$ the fan corresponding to $X_4^0$. Then $\G(\Sigma_3^0)=\{x_0,x_+,x_-,y_+,y_-,z_1\}$, while $\G(\Sigma_4^0)=\{x_0,x_+,x_-,y_+,y_-,z_1,z_2\}$. Especially, $X_3^0$ is a $F_1$-bundle over ${\bf P}^1$, while $X_4^0$ is a $S_7$-bundle over ${\bf P}^1$. Moreover, we have
$$\left(-K_{X_3^0}\right)^3=52\mbox{ and }\left(-K_{X_4^0}\right)^3=38.$$

There exists a sequence of equivariant blow-ups along curves
$$X_4^0\longrightarrow X_3^0\longrightarrow{\bf P}_{{\bf P}^1}\left(\mathcal{O}_{{\bf P}^1}\oplus\mathcal{O}_{{\bf P}^1}\oplus\mathcal{O}_{{\bf P}^1}(2)\right).$$
On the other hand, there exists a sequence of blow-ups along curves
$$Y_4^0\longrightarrow Y_3^0\longrightarrow{\bf P}_{{\bf P}^1}\left(\mathcal{O}_{{\bf P}^1}\oplus\mathcal{O}_{{\bf P}^1}(1)\oplus\mathcal{O}_{{\bf P}^1}(1)\right),$$
where $Y_3^0$ is a toric Fano $3$-fold of Picard number $3$ and of type no.31 on the table in Mori-Mukai \cite{mori1}, while $Y_4^0$ is a Fano $3$-fold, which is not toric, of Picard number $4$ and of type no.8 on the table in Mori-Mukai \cite{mori1}. $X_3^0$ and $X_4^0$ are deformed to $Y_3^0$ and $Y_4^0$ under small deformations, respectively. Moreover, ${\bf P}_{{\bf P}^1}\left(\mathcal{O}_{{\bf P}^1}\oplus\mathcal{O}_{{\bf P}^1}\oplus\mathcal{O}_{{\bf P}^1}(2)\right)$ is deformed to ${\bf P}_{{\bf P}^1}\left(\mathcal{O}_{{\bf P}^1}\oplus\mathcal{O}_{{\bf P}^1}(1)\oplus\mathcal{O}_{{\bf P}^1}(1)\right)$, though ${\bf P}_{{\bf P}^1}\left(\mathcal{O}_{{\bf P}^1}\oplus\mathcal{O}_{{\bf P}^1}(1)\oplus\mathcal{O}_{{\bf P}^1}(1)\right)$ is not Fano (see Ashikaga-Konno \cite{ashikaga1}, Harris \cite{harris1} and Nakamura \cite{nakamura1}).

\begin{Cor}\label{p1p1}

Let $X$ be a toric weakened Fano $3$-fold. If there exists a $(0,2)$-type contraction whose exceptional divisor is isomorphic to ${\bf P}^1\times{\bf P}^1$, then the exceptional divisors of other $(0,2)$-type contractions are also isomorphic to ${\bf P}^1\times{\bf P}^1$.

\end{Cor}

\bigskip

$({\rm II})$ $E\cong F_1\ (a=1).$

\bigskip

In this case, we have to consider the following (see Remark \ref{f1only}).

\begin{Lem}\label{korochan}

$\{y_+,y_-\}$ is a primitive collection. Moreover, $z:=y_++y_-$ is contained in $\G(\Sigma)$ and the primitive relation of $\{y_+,y_-\}$ is $y_++y_-=z$.

\end{Lem}

To prove Lemma \ref{korochan}, we need the following.

\begin{Lem}\label{tete}

Let $X$ be a toric weakened Fano $3$-fold and $\Sigma$ the corresponding fan. Any primitive collection $P=\{x_1,x_2\}$ such that its primitive relation is $x_1+x_2=2y$ $(y\in\G(\Sigma))$ and $\{x_1,y\}$ is a ${\bf Z}$-basis of $N$ is extremal.

\end{Lem}

\proof
If $P$ is not extremal, by the same argument as the proof of Lemma \ref{key}, there exists an extremal primitive relation $x'+x''=2y$ whose corresponding contraction is $(0,2)$-type, where $x',x''\in\G(\Sigma)$. By the assumption $\{x_1,y\}$ is a ${\bf Z}$-basis of $N$, ${\bf R}_{\geq 0}x_1+{\bf R}_{\geq 0}x_2$ does not contain ${\bf R}_{\geq 0}x'+{\bf R}_{\geq 0}x''$. This is impossible, because $X$ is a weak Fano variety and $x'+x''=2y$ corresponds to a $(0,2)$-type contraction.
\qed

\bigskip

Proof of Lemma \ref{korochan}. \quad
Suppose $\{y_+,y_-\}$ is not a primitive collection. Then we have two primitive relations
$$x_-+y_++y_-=x_0\mbox{ and }x_0+y_++y_-=x_+.$$
The completeness of $X$, there exists
$$z=\pmatrix{\alpha\cr\beta\cr\gamma}\in I$$
such that $\alpha<0$. Obviously, $\{x_+,z\}$ is a primitive collection. If $I\subset\mathcal{S}$ then by the same argument as in the proof of Lemma \ref{wataru}, the corresponding primitive relation have to be $x_++z=0$. So we have
$$x_-+z=2\pmatrix{0\cr -1\cr 0},\ \pmatrix{0\cr -1\cr 0}\in\G(\Sigma).$$
By Lemma \ref{tete}, this primitive collection is extremal. However, the corresponding contraction is not of $(0,2)$-type. So suppose $I\subset\mathcal{T}$. This case is also impossible by the same argument as in the proof of Lemma \ref{wataru}. Therefore $\{y_+,y_-\}$ is a primitive collection.
\qed

\bigskip

Since the Picard number of $X$ is not three, we have $I\subset\mathcal{S}$ and $X$ is a toric surface bundle over ${\bf P}^1$ by $(1)$ in Corollary \ref{bdl} and Lemma \ref{korochan}. The fan corresponding to a fiber of this bundle structure is in $\mathcal{S}$. By the classification of toric weak del Pezzo surfaces in Section \ref{surface}, in this case, we get two new toric weakened Fano $3$-folds $X_4^1$ and $X_5^1$ whose Picard numbers are $4$ and $5$, respectively. The other cases are impossible since there exists a crepant contraction which is not of $(0,2)$-type. Put
$$z_1=\pmatrix{0\cr 1\cr 0},\ z_2=\pmatrix{-1\cr 0\cr 0},\ z_3=\pmatrix{0\cr -1\cr 0},$$
$\Sigma_4^1$ the fan corresponding to $X_4^1$ and $\Sigma_5^1$ the fan corresponding to $X_5^1$. Then $\G(\Sigma_4^1)=\{x_0,x_+,x_-,y_+,y_-,z_1,z_2\}$, while $\G(\Sigma_5^1)=\{x_0,x_+,x_-,y_+,y_-,z_1,z_2,z_3\}$. Especially, $X_4^1$ is a $W_3$-bundle over ${\bf P}^1$, while $X_5^1$ is a $W^1_4$-bundle over ${\bf P}^1$ (see Table 1 in Section \ref{surface}). Moreover, we have
$$\left(-K_{X_4^1}\right)^3=46\mbox{ and }\left(-K_{X_5^1}\right)^3=36.$$

There exists an equivariant blow-up along a curve
$$X_5^1\longrightarrow X_4^1.$$
On the other hand, there exists a blow-up along a curve
$$Y_5^1\longrightarrow Y_4^1,$$
where $Y_4^1$ is a toric Fano $3$-fold of Picard number $4$ and of type no.12 on the table in Mori-Mukai \cite{mori1}, while $Y_5^1$ is a Fano $3$-fold, which is not toric, of Picard number $5$ and of type no.2 on the table in Mori-Mukai \cite{mori1}. $X_4^1$ and $X_5^1$ are deformed to $Y_4^1$ and $Y_5^1$ under small deformations, respectively.

\begin{Cor}\label{f1f1}

Let $X$ be a toric weakened Fano $3$-fold. If there exists a $(0,2)$-type contraction whose exceptional divisor is isomorphic to $F_1$, then $X$ is isomorphic to either $X_4^1$ or $X_5^1$.

\end{Cor}

\bigskip


$({\rm III})$ $E\cong F_2\ (a=2).$

\bigskip

We have a primitive relation $y_++y_-=2z$, where
$$z=\pmatrix{0\cr 1\cr 0}\in\G(\Sigma).$$
So the Picard number of $X$ is not three. Therefore, $I\subset\mathcal{S}$ and $X$ is a toric surface bundle over ${\bf P}^1$ by $(1)$ in Corollary \ref{bdl}. The fan corresponding to a fiber of this bundle structure is in $\mathcal{S}$. On the other hand, by Lemma \ref{tete}, the primitive relation $y_++y_-=2z$ is extremal, therefore, the corresponding contraction is of $(0,2)$-type. Moreover, its exceptional divisor have to be isomorphic to $F_2$ by Corollaries \ref{p1p1} and \ref{f1f1}. This is impossible. Thus there exists no $(0,2)$-type contraction whose exceptional divisor is isomorphic to $F_2$ on any toric weakened Fano $3$-fold.

\begin{Rem}

{\rm
There is an example of a general weakened Fano $3$-fold which has a $(0,2)$-type contraction whose exceptional divisor is isomorphic to $F_2$ (see Minagawa \cite{minagawa1}).
}

\end{Rem}

Thus, we get the classification of toric weakened Fano $3$-folds.

\begin{Thm}\label{class}

There exist exactly $15$ smooth toric weakened Fano $3$-folds up to isomorphisms. There are following three cases$:$
\begin{enumerate}
\item ${\bf P}^1\times X'$, where $X'$ is a toric weak del Pezzo surface but not a del Pezzo surface, that is, toric weakened del Pezzo surface$:$ ${\bf P}^1\times F_2$, ${\bf P}^1\times W_3$, ${\bf P}^1\times W^1_4$, ${\bf P}^1\times W^2_4$, ${\bf P}^1\times W^3_4$, ${\bf P}^1\times W^1_5$, ${\bf P}^1\times W^2_5$, ${\bf P}^1\times W^1_6$, ${\bf P}^1\times W^2_6$, ${\bf P}^1\times W^3_6$ and ${\bf P}^1\times W_7$.
\item Toric del Pezzo surface bundles over ${\bf P}^1:$ $X^0_3$ and $X^0_4$.
\item Toric weakened del Pezzo surface bundles over ${\bf P}^1$ which are not decomposed into a direct product of ${\bf P}^1$ and a toric weakened del Pezzo surface$:$ $X^1_4$ and $X^1_5$.
\end{enumerate}

\end{Thm}

\bigskip

\begin{flushleft}
\begin{sc}
Department of Mathematics, \\
Tokyo Institute of Technology, \\
Oh-Okayama, Meguro, Tokyo, Japan
\end{sc}

\medskip
{\it E-mail address}: $\mathtt{hirosato@math.titech.ac.jp}$
\end{flushleft}


\begin{thebibliography}{99}
\bibitem{ashikaga1} T. Ashikaga and K. Konno, Algebraic surfaces of general type with $c_1^{2}=3p_{g}-7$, Tohoku Math. J. 42 (1990), 517--536.

\bibitem{batyrev3} V. V. Batyrev, On the classification of smooth projective toric varieties, Tohoku Math. J. 43 (1991), 569--585.

\bibitem{batyrev4} V. V. Batyrev, On the classification of toric Fano 4-folds, Algebraic geometry, 9, J. Math. Sci. (New York) 94 (1999), 1021--1050.

\bibitem{fulton1} W. Fulton, Introduction to Toric Varieties, Ann. of Math. Studies 131, Princeton Univ. Press, 1993.

\bibitem{harris1} J. Harris, A bound on the geometric genus of projective varieties, Ann. Sci. Norm. Sup. Pisa, IV. Ser VIII (1981), 35--68.

\bibitem{koelman1} R. Koelman, The number of moduli of families of curves on toric surfaces, Thesis, Univ. Nijmegen, 1991.


\bibitem{minagawa1} T. Minagawa, On classification of weakened Fano $3$-folds with $B_2(X)=2$, in Proc. of algebraic geometry symposium (Kinosaki, Oct. 2000).

\bibitem{minagawa2} T. Minagawa, Global smoothing of singular weak Fano $3$-folds, preprint (1999).

\bibitem{mori1} S. Mori and S. Mukai, Classification of Fano $3$-folds with $B_2\geq 2$, Manuscripta Math. (1981), no. 36, 147--162.

\bibitem{nakamura1} I. Nakamura, Global deformations of ${\bf P}^{2}$-bundles over ${\bf P}^{1}$, J. Math. Kyoto Univ. 38 (1998), 29--54.

\bibitem{oda2} T. Oda, Convex Bodies and Algebraic Geometry---An introduction to the theory of toric varieties, Ergeb. Math. Grenzgeb. (3), Vol. 15, Springer-Verlag, Berlin, Heidelberg, New York, London, Paris, Tokyo, 1988.

\bibitem{reid1} M. Reid, Decomposition of toric morphisms, in Arithmetic and Geometry, papers dedicated to I. R. Shafarevich on the occasion of his 60th birthday (M. Artin and J. Tate, eds.), vol. II, Geometry, Progress in Math. 36, Birkh\"{a}user, Boston, Basel, Stuttgart, 1983, 395--418.

\bibitem{sato1} H. Sato, Toward the classification of higher-dimensional toric Fano varieties, Tohoku Math. J. 52 (2000), 383--413.


\end{thebibliography}
\end{document}